\documentclass [12 pt, twoside]{amsart}
\usepackage{enumerate}

 \usepackage{amsmath,amsfonts,amsthm,amssymb}
 \usepackage[all]{xy}

\newtheorem{thm}{Theorem}[section]
\newtheorem{prop}{Proposition}[section] 
\newtheorem{lem}{Lemma}[section] 
\newtheorem{cor}{Corollary}[section]

\newtheorem{que}{Question}[section]

\newcommand{\bbb}[1]{\mbox{\boldmath$#1$}}

\theoremstyle{definition}
\newtheorem*{Proof}{Proof}

\newcommand{\ga} {{\gamma}}

\newcommand{\ld} {{\ldots}}
\newcommand{\sm} {{\smallsetminus}}
\newcommand{\thi} {{\theta}}
\newcommand{\de} {{\delta}}
\newcommand{\De} {{\varDelta}}
\newcommand{\si} {{\sigma}}
\newcommand{\Si} {{\varSigma}}
\newcommand{\la} {{\lambda}}
\newcommand{\La} {{\varLambda}}

\newcommand{\e} {{\varepsilon}}

\newcommand{\mi} {{\mu}}

\newcommand{\dis}{\displaystyle}
\newcommand{\ssum}{\sum\limits}
\newcommand{\dlim}{\displaystyle\lim}

\newcommand{\cp}{{\mathcal{P}}}

\newcommand{\cb}{{\mathcal{B}}}

\newcommand{\ra}{{\rightarrow}}
\newcommand{\fa}{{\forall}}

\newcommand{\qb}{$\quad\blacksquare$}
\def\1{\it1\hspace*{-0.150cm}{\footnotesize{I}}}

\def\C{{\mathbb{C}}}
\def\Q{{\mathbb{Q}}}
\def\N{{\mathbb{N}}}

\numberwithin{equation}{section}

\begin{document}

\title[ Common hypercyclic functions ]{ Common hypercyclic functions for translation operators with large gaps}

\author[Nikos Tsirivas]{Nikos Tsirivas}
\address{Department of Mathematics and Applied Mathematics, University of Crete, GR-700 13 Heraklion, Crete, Greece}
\email{tsirivas@uoc.gr}
\thanks{The research project is implemented within the framework of the Action "Supporting Postdoctoral Researchers" of
the Operational Program "Education and Lifelong Learning" (Action's Beneficiary: General Secretariat for
Research and Technology), and is co-financed by the European Social Fund (ESF) and the Greek State.}


\date{}


\begin{abstract}
Let $H(\C)$ be the set of entire functions endowed with the topology of local uniform convergence. Fix a
sequence of non-zero complex numbers $(\la_n)$, $|\la_n|\to +\infty$, which satisfies the following property:for
every $M>0$ there exists a subsequence $(\mi_n)$ of $(\la_n)$ such that \smallskip

(i) $|\mi_{n+1}|-|\mi_n|>M$ for every $n=1,2,\ld$ and  \smallskip

(ii) $\ssum^{+\infty}_{n=1}\dfrac{1}{|\mi_n|}=+\infty$ .\smallskip

We prove that there exists a residual set $G\subset H(\C)$ such that for every $f\in G$ and every non-zero
complex number $a$ the set $\{ f(z+\la_na):n=1,2,\ldots \}$ is dense in $H(\C)$. This answers in the affirmative
Question 1 in \cite{7} and it also provides an extension of a theorem due to Costakis and Sambarino in \cite{6}.
\end{abstract}

\maketitle

\section{Introduction}
\noindent We start by fixing some standard notation and terminology. The symbols $\N=$ $\{1,2,\ld\}$, $\Q$,
$\mathbb{R}$, $\mathbb{C}$ stand for the sets of natural, rational, real and complex numbers respectively. By
$H(\C)$ we denote the set of entire functions endowed with the topology of local uniform convergence. For a
subset $A$ of $H(\C)$, $\overline{A}$ denotes the closure of $A$ with respect to the topology of local uniform
convergence. Let $X$ be a topological vector space. A subset $G$ of a $X$ is called $G_{\de}$ if it can be
written as a countable intersection of open sets in $X$ and a subset $Y$ of $X$ is called residual if it
contains a $G_{\de}$ and dense subset of $X$. The symbol $\infty$ whenever appears in the present work denotes
the complex infinity.

Let $(T_n:X\ra X)$ be a sequence of continuous linear operators on a topological vector space $X$. If
$(T_n(x))_{n\ge1}$ is dense in $X$ for some $x\in X$, then $x$ is called \textit{hypercyclic} for $(T_n)$ and we
say that $(T_n)$ is \textit{hypercyclic} \cite{2}, \cite{10}. The symbol $HC(\{ T_n\})$ stands for the
collection of all hypercyclic vectors for $(T_n)$. In the case where the sequence $(T_n)$ comes from the
iterates of a single operator $T:X\to X$, i.e. $T_n:=T^n$, then we simply say that $T$ is \textit{hypercyclic}
and $x$ is \textit{hypercyclic} for $T$. If $T:X\to X$ is hypercyclic then the symbol $HC(T)$ stands for the
collection of all hypercyclic vectors for $T$. A simple consequence of Baire's category theorem is that for
every continuous linear operator $T$ on a separable topological vector space $X$, if $HC(T)$ is non-empty then
it is necessarily ($G_{\de}$ and) dense. For an account of results on the subject of hypercyclicity we refer to
the recent books \cite{2}, \cite{10}, see also the very influential survey article \cite{Grosse2}.

In the present work we deal with translation operators. For every $a\in \mathbb{C} \setminus \{ 0\}$ consider
the translation operator $T_a:H(\C) \to H(\C)$ defined by $$T_a(f)(z)=f(z+a), \,\,\, f\in H(\C).$$ An old result
of Birkhoff \cite{4} says that there exist entire functions the integer translates of which are dense in the
space of all entire functions endowed with the topology of local uniform convergence. In other words $T_1$ is
hypercyclic. Actually, it is not difficult to see that for every $a\in \mathbb{C} \setminus \{ 0\}$, $T_a$ is
hypercyclic and hence $HC(T_a)$ is $G_{\de}$ and dense in $H(\C)$. Costakis and Sambarino \cite{6} strengthened
Birkhoff's result by showing that the family $\{ T_a| \,\, a\in \mathbb{C} \setminus \{ 0\} \}$ has a residual
set of common hypercyclic vectors i.e., the set $\bigcap_{a\in \mathbb{C}\setminus \{ 0\}} HC(\{ T_{na}\} ) $is
residual in $H(\C)$. In particular, it is non-empty. Of course, what makes their result non-trivial is the
uncountable range of $a$. At this point, let us mention a relevant observation due to Bayart and Matheron,
\cite{2}, \cite{3}: suppose $X$ is a Fr\'{e}chet space and $\{ S_{a,n} |\, a\in A, n\in \mathbb{N} \}$ is a
collection of sequences of continuous linear operators on $X$, labelled by the elements $a$ of a set $A$. If $A$
is a $\sigma$-compact topological space, the maps $a \to S_{a,n}$ are $SOT$-continuous and each sequence
$(S_{a,n})_{n\in \mathbb{N}}$ has a dense set of hypercyclic vectors then either $\bigcap_{a\in A}HC(\{ S_{a,n}
\})=\emptyset$ or $\bigcap_{a\in A}HC(\{ S_{a,n} \})$ is a dense $G_{\delta}$-set in $X$. This observation
applies to all the collections of operators considered in our work.

Let us now come to the main subject of our paper. Recall that the set $\bigcap_{a\in \mathbb{C}\setminus \{ 0\}}
HC(\{ T_{na}\} ) $ is residual in $H(\C)$, \cite{6}. Subsequently, Costakis \cite{7} asked whether, in this
result, the sequence $(n)$ can be replaced by more general sequences $(\la_n)$ of non-zero complex numbers. In
this direction Costakis \cite{7} showed that, if the sequence $(\la_n)$ satisfies the following condition
$(\Si)$: for every $M>0$ there exists a subsequence $(\mi_n)$ of $(\la_n)$ such that
\smallskip

(i) $|\mi_{n+1}|-|\mi_n|>M$ for every $n=1,2,\ld$ and  \smallskip

(ii) $\ssum^{+\infty}_{n=1}\dfrac{1}{|\mi_n|}=+\infty$ ,\smallskip

then the desired conclusion holds if we restrict attention to $a\in C(0,1):=\{z\in\C/|z|=1\}$, that is the set
$\bigcap_{a\in C(0,1)} HC(\{ T_{\la_na} \})$ is residual in $H(\C)$. In view of the above, Costakis led to the
following question, see Question 1 in \cite{7}.

\begin{que} \label{Q1}
Let $(\la_n)$ be a sequence of non-zero complex numbers tending to infinity which also satisfies condition
$(\Si)$. Is it true that the set $\bigcap_{a\in \mathbb{C}\setminus \{ 0\}} HC(\{ T_{\la_na} \})$ is residual in
$H(\C)$, hence non-empty?
\end{que}

Our main task is to give an affirmative answer to Question \ref{Q1} by proving the following
\begin{thm}\label{thm2.1}
Fix a sequence of non-zero complex numbers $\La=(\la_n)$ that tends to infinity and satisfies the above
condition $(\Si)$.\ Then $\bigcap\limits_{a\in\C\sm\{0\}}HC(\{T_{\la_na}\})$ is a $G_\delta$ and dense subset of
$H(\C)$.
\end{thm}
It is worth to mention here that one is forced to impose certain natural restrictions on the sequence $(\la_n)$
in order to conclude that the set $\bigcap\limits_{a\in\C\sm\{0\}}HC(\{T_{\la_na}\})$ is non-empty. Indeed, in
\cite{8} the authors show that if $\liminf_n\frac{|\la_{n+1}|}{|\la_n|}>2$ then
$\bigcap\limits_{a\in\C\sm\{0\}}HC(\{T_{\la_na}\})=\emptyset $. In particular,
$\bigcap\limits_{a\in\C\sm\{0\}}HC(\{T_{e^na}\})=\emptyset $. However, for sequences $(\lambda_n)$ with
$$1<\liminf_n\frac{|\la_{n+1}|}{|\la_n|}\leq 2$$ it is not known whether
$$\bigcap\limits_{a\in\C\sm\{0\}}HC(\{T_{\la_na}\})=\emptyset, $$ although it is plausible to conjecture that this
is the case. In particular, we do not know what happens when $\lambda_n =2^n$ or $\lambda_n=(3/2)^n$. This work
can be seen as a try to understand the nature of this restriction. In any case, it seems a quite difficult
problem to fully characterize the sequences $(\la_n)$ for which the conclusion of Theorem \ref{thm2.1} holds.

We stress that Theorem \ref{thm2.1} complements the main result from our recent work in \cite{Tsi}. In
\cite{Tsi} we showed that the conclusion of Theorem \ref{thm2.1} holds for sequences $(\la_n)$ satisfying
another type of condition different from $(\Si)$; this condition, which we call it $(\Si ')$, is also not very
restrictive, in the sense that it still allows sequences $(\la_n)$ with ``large gaps". To avoid extra notation
and to keep the introduction in a compact form, we postpone the definition of condition $(\Si ')$ till section
6. We note that although sequences of polynomial type of degree bigger than one, such as $(n^2)$, $(n^3)$,
$(n+n^3)$, $(n^4+n^5)$ and so on, clearly do not satisfy condition $(\Si )$ they do satisfy $(\Si ')$. On the
other hand there exist sequences satisfying $(\Si ')$ which do not satisfy $(\Si )$. However, there exist
sequences satisfying both conditions $(\Si )$ and $(\Si ')$. All these are explained in full detail in Section
6.

A few words about the proof of Theorem \ref{thm2.1}. Of course the main argument uses Baire's category theorem,
but in order to do so the first and most difficult thing is to construct a suitable two dimensional partition on
a given sector of the plane. After, to each point of the partition we assign a suitable closed disk of constant
radius so that these disks are pairwise disjoint and their union almost fills the sector. Having done these
steps we are ready for the final argument which involves a standard use of Runge's or Mergelyan's approximation
theorem along with Baire's theorem. It is important to say that in our framework one cannot use Ansari's theorem
\cite{Ans}, as Costakis and Sambarino did in their proof, since now the sequence $(\la_n)$ lacks the semigroup
structure, i.e. $\la_n +\la_m \neq \la_{n+m}$ in general. Actually, this was the reason that led us to seek
higher order partitions in order to make things work. Overall, we elaborate on the work of Costakis and
Sambarino and we offer a general strategy how to construct two dimensional partitions relevant to our problem.
In general, our proof shares certain similarities with the proof of the main result in \cite{Tsi} and so we feel
that the interested reader will get a more clear and integrated picture by reading in parallel the present paper
and paper \cite{Tsi}. However, the methods of constructing the partitions in the present paper and \cite{Tsi}
differentiate drastically. The reason for this, is that always the partition reflects the structure of the
sequence $(\lambda_n)$. The construction of the partition in \cite{Tsi} is very tight and quite delicate and
comes from our effort to deal firstly with the most natural sequence which fails condition $(\Si)$, namely the
sequence $(n^2)$. It is also evident that there is a huge distance between sequences satisfying condition $(\Si
)$ and the sequences satisfying condition $(\Si ')$, see section 6. Of course, it would be desirable to exhibit
a condition and a corresponding partition, if any, which imply the main result of the present paper as well as
the main result in \cite{Tsi}. Unfortunately, this is unclear to us.

There are several recent results concerning either the existence or the non-existence of common hypercyclic
vectors for uncountable families of operators, such as weighted shifts, adjoints of multiplication operators,
differentiation and composition operators; see for instance, \cite{AbGo}, \cite{Ba1}-\cite{BaGrMo}, \cite{Ber}
\cite{ChaSa1}-\cite{ChaSa5}, \cite{5}-\cite{8}, \cite{GaPa}, \cite{10}, \cite{LeMu}, \cite{Math}, \cite{San},
\cite{Shka}, \cite{13}, \cite{Tsi}.

Our paper is organized as follows. The proof of Theorem \ref{thm2.1} has several steps and occupies Sections
2-5. Finally, in Section 6 we compare Theorem \ref{thm2.1} with the main result from \cite{Tsi} and we exhibit
examples of sequences which illustrate our main theorem.

\section{A reduction of Theorem \ref{thm2.1}}\label{sec2}
\noindent

Let us now describe the steps for the proof of Theorem \ref{thm2.1}. Consider the sectors
\[
S^k_n:=\bigg\{a\in\C|\exists\;r\in\bigg[\frac{1}{n},n\bigg]\;\text{and}\;t\in\bigg[\frac{k}{4},
\frac{k+1}{4}\bigg]\;\text{such that}\;a=re^{2\pi it}\bigg\}
\]
for $k=0,1,2,3$ and $n=2,3,\ld\,$. Since
\[
\bigcap_{a\in\C\sm\{0\}}HC(\{ T_{\la_na} \} )=\bigcap^3_{k=0}\bigcap^{+\infty}_{n=2} \bigcap_{a\in S^k_n}HC(\{
T_{\la_na}\} ),
\]
an appeal of Baire's category theorem reduces Theorem \ref{thm2.1} to the following.
\begin{prop}\label{prop2.1}
Fix a sequence $(\lambda_n)$ of non-zero complex numbers that tends to infinity which satisfies the above condition $(\Si)$. Fix
four real numbers $r_0,R_0,\thi_0,\thi_T$ such that $0<r_0<1<R_0<+\infty$, $0\le\thi_0<\thi_T \leq 1$,
$\thi_T-\thi_0=\dfrac{1}{4}$ and consider the sector $S$ defined by
\[
S:=\{a\in\C | \,\,\text{there exist} \,\,r\in[r_0,R_0] \ \ \text{and} \ \ t\in[\thi_0,\thi_T] \ \ \text{such
that} \ \ a=re^{2\pi it}\}.
\]
Then $\bigcap\limits_{a\in S}HC(\{ T_{\la_na}\} )$ is a $G_\delta$ and dense subset of $H(\mathbb{C})$.
\end{prop}
For the proof of Proposition \ref{prop2.1} we introduce some notation which will be carried out throughout this
paper. Let $(p_j)$, $j=1,2,\ld$ be a dense sequence of $H(\mathbb{C})$, (for instance, all the polynomials in
one complex variable with coefficients in $\Q+i\Q$). For every $m,j,s,k\in\N$ we consider the set
$$E(m,j,s,k):= \!\Big\{f\!\in\!H(\mathbb{C})\,\, | \forall a\!\in\! S \,\,\,\exists  n\!\in\!\N, n\!\le\! m :
\dis\sup_{|z|\le k}|f(z\!+\!\la_na)\!-\!p_j(z)|\!<\!\frac{1}{s}\Big\}.$$
By Baire's category theorem and the three lemmas stated below, Proposition \ref{prop2.1} readily follows.
\begin{lem}\label{lem2.1}
\[
\bigcap_{a\in S}HC(\{ T_{\la_na}\} )=\bigcap^{+\infty}_{j=1}\bigcap^{+\infty}_{s=1}\bigcap^{+\infty}_{k=1}
\bigcup^{+\infty}_{m=1}E(m,j,s,k).
\]
\end{lem}
\begin{lem}\label{lem2.2}
For every $m,\!j,\!s,\!k\!\in\!\N$ the set $E(m,\!j,\!s,\!k)$ is open in $H(\mathbb{C})$.
\end{lem}
\begin{lem}\label{lem2.3}
For every $j,\!s,\!k\in\N$ the set $\bigcup\limits^{+\infty}_{m=1}E(m,j,s,k)$ is dense in $H(\mathbb{C})$.
\end{lem}

The proof of Lemma \ref{lem2.1} is in \cite{Tsi}. The proof of Lemma \ref{lem2.2} is similar to that in Lemma 9
of \cite{6} and it is omitted.

We now move on to Lemma \ref{lem2.3}. This lemma is the heart of our argument and its proof occupies the next
three sections.

\section{Construction of the partition of the sector $\bbb{S}$}\label{sec3}
\noindent

For the sequel we fix four positive numbers $c_1,c_2,c_3,c_4$ such that $c_1>1$, $c_2\in(0,1)$, $c_3>1$,
$c_4>1$, where $c_3:=\dfrac{c_4}{r_0c_2}$, $c_1:=4(c_3+1)$. We also consider four positive real numbers
$\thi_0,\thi_T,r_0,R_0$ as in Proposition \ref{prop2.1} and a sequence $\La=(\la_n)$ of non zero complex numbers
which satisfies condition $(\Si)$ and such that $\la_n\ra\infty$ as $n\ra+\infty$. After the definition of the
above numbers we fix a subsequence $(\mi_n)$ of $(\la_n)$ such that:
\[
|\mi_n|>c_1, |\mi_{n+1}|-|\mi_n|>c_1 \ \ \text{for every} \ \ n=1,2,\ld \ \ \text{and} \ \
\sum^{+\infty}_{k=1}\frac{1}{|\mi_k|}=+\infty.
\]
\subsection{Step 1. Partitions of the inverval $\bbb{[\thi_0,\thi_T]}$} \label{subsec3.1}
\noindent

In this step we succeed the elementary structure of our construction. The following two steps are based in this
first one. For every positive integer $m$ we shall construct a corresponding partition $\De_m$ of
$[\thi_0,\thi_T]$. So, let $m\in\N$ be fixed.\\

The condition $\ssum^{+\infty}_{n=1}\dfrac{1}{|\mi_n|}=+\infty$ implies that for every positive integer $m=1,2,\ld$
there exists the minimum natural number $m_1(m)$ such that:
\begin{eqnarray}
\sum^{m_1(m)}_{k=m}\frac{1}{|\mi_k|}>c_3\cdot\frac{1}{|\mi_m|}.  \label{eq3.1}
\end{eqnarray}
Clearly $m_1(m)\ge m+1$ for every $m=1,2,\ld$ because $c_3>1$. We define the numbers $\thi_0^{(m)}:=\thi_0$,
$\thi^{(m)}_1:=\thi^{(m)}_0+\dfrac{c_2}{|\mi_m|}$, $\thi^{(m)}_2:=\thi^{(m)}_1+\dfrac{c_2}{|\mi_{m+1}|},\ld$,
$\thi^{(m)}_{m_1(m)-m+1}:=\thi^{(m)}_{m_1(m)-m}+\dfrac{c_2}{|\mi_{m_1(m)}|}$, or generally:
\begin{eqnarray}
\thi^{(m)}_{n+1}:=\thi^{(m)}_n+\frac{c_2}{|\mi_{m+n}|}, \ \ n=0,1,\ldots , m_1(m)-m,  \label{eq3.2}
\end{eqnarray}
where $m_1(m)-m\ge1$. Define $$\si_m:=\thi^{(m)}_{m_1(m)-m+1}-\thi_0.$$ Now let any positive integer $\nu$ with
$$\nu>m_1(m)-m+1.$$
For such a $\nu$ there exists a unique pair $(k,j)\in\N^2$, where $j\in\{0,1,\ld,m_1(m)-m\}$, such that:
\[
\nu=k(m_1(m)-m+1)+j.
\]
We define
\[
\thi^{(m)}_\nu:=\thi^{(m)}_j+k\si_m.
\]
It is obvious that $\dlim_{\nu\ra+\infty}\thi^{(m)}_\nu=+\infty$ and the sequence $(\thi^{(m)}_\nu)_\nu$ is
strictly increasing, in respect to $\nu$. So there exists a maximum natural number $\nu_m\in\N$ such that
$\thi^{(m)}_{\nu_m}\le\thi_T$. We set $$\De_m:=\{\thi^{(m)}_0,\thi^{(m)}_1,\ld,\thi^{(m)}_{\nu_m}\}.$$ It holds
that $\nu_m\ge m_1(m)-m+1$ ( see Lemma \ref{lemma3.4.1}).
\subsection{Step 2. Partitions of the arc $\bbb{\phi_r([\theta_0 ,\theta_T])}$} \label{subsec3.2}
\noindent

Consider the function $\phi :[\theta_0 ,\theta_T ]\times (0 ,+\infty )\to \mathbb{C}$ given by
$$\phi (t,r):=re^{2\pi it},\,\,\, (t,r)\in [\theta_0 ,\theta_T ]\times (0 ,+\infty ) $$ and for every $r> 0$ we
define the corresponding curve $\phi_r:[\theta_0 ,\theta_T ] \to \mathbb{C}$ by
$$\phi_r(t):=\phi (t,r), \,\,\, t\in [\theta_0 ,\theta_T ].$$
For any given positive integer $m$, $\phi_r(\Delta_m)$ is a partition of the arc $\phi_r([\theta_0 ,\theta_T])$,
where $\Delta_m$ is the partition of the interval $[\theta_0 ,\theta_T]$ constructed in Step 1. For every $r>0$,
$m\in \mathbb{N}$ define
$$\cp_0^{r,m}:=\phi_r(\Delta_m)$$
which we call partition of the arc $\phi_r([\theta_0 ,\theta_T])$ with height $r$, density $m$ and order $0$.
\subsection{Step 3. The final partition}  \label{subsec3.3}
\noindent

Consider the partition $\cp^{r_0,1}_0$ from the previous step, Step 2 and set
\begin{eqnarray}
r_1:=r_0+\frac{c_2}{|\mi_{m_1(1)}|}.  \label{eq3.3}
\end{eqnarray}
After, we consider the partition $\cp_0^{r_1,m_1(1)+1}$ and we set $$m_2:=m_1(m_1(1)+1),$$
$$r_2:=r_1+\dfrac{c_2}{|\mi_{m_2}|}.$$ Inductively we define two sequences $(r_\nu)$, $\nu=0,1,2,\ld$, $(m_\nu)$,
$\nu=2,\ld$, as follows: $r_0,r_1,r_2$ and $m_2$ are as above, see (\ref{eq3.3}). Suppose that we have
constructed the numbers $m_\nu$, $r_\nu$ for some $\nu\ge2$. Then, taking into account the partition
$\cp_0^{r_\nu,m_\nu+1}$, we set
\begin{eqnarray}
m_{\nu+1}=m_1(m_\nu+1) \label{eq3.4}
\end{eqnarray}
and
\begin{eqnarray}
r_{\nu+1}:=r_\nu+\frac{c_2}{|\mi_{m_{\nu+1}}|}.  \label{eq3.5}
\end{eqnarray}
For the next step, consider the partition $\cp_0^{r_{\nu+1},m_{\nu+1}+1}$. We will prove in the next subsection
that $\dlim_{\nu\ra+\infty}r_\nu=+\infty$. Therefore there exists a maximum natural number $\nu_0\in\N$ such
that $r_{\nu_0}\le R_0$ because the sequence $(r_\nu)$ is strictly increasing. In view of the above, we define
\[
\cp:=\cp_0^{r_0,1}\cup\bigg(\bigcup^{\nu_0}_{\nu=1}\cp_0^{r_\nu,m_\nu+1}\bigg),
\]
which is the desired partition of our sector $S$.
\subsection{Properties of the partitions} \label{subsec3.4}
\begin{lem} \label{lemma3.4.1}
Let some fixed $m\in\N$. Then
\[
\si_m=\thi^{(m)}_{m_1(m)-m+1}-\thi_0<\frac{1}{4}.
\]
In particular, $\nu_m\ge m_1(m)-m+1$.
\end{lem}
\begin{Proof}
By the definition of the numbers $\thi^m_j$, $j=0,1,\ld,m_1(m)-m+1$ we have
\begin{eqnarray}
\thi^{(m)}_{m_1(m)-m+1}-\thi_0=c_2\cdot\sum^{m_1(m)}_{k=m}\frac{1}{|\mi_k|},  \label{eq3.6}
\end{eqnarray}
and by the definition of the number $m_1(m)$ it follows that
\begin{eqnarray}
\sum^{m_1(m)}_{k=m}\frac{1}{|\mi_k|}\le
c_3\cdot\frac{1}{|\mi_m|}+\frac{1}{|\mi_{m_1(m)}|}<(c_3+1)\frac{1}{|\mi_m|}.  \label{eq3.7}
\end{eqnarray}
Our hypotheses imply $c_1=4(c_3+1)$ and $|\mi_m|>c_1=4(c_3+1)>4c_2(c_3+1)$, because $c_2\in(0,1)$. This gives
\begin{eqnarray}
\frac{c_3+1}{|\mi_m|}<\frac{1}{4c_2}.  \label{eq3.8}
\end{eqnarray}
Thus, (\ref{eq3.6}), (\ref{eq3.7}) and (\ref{eq3.8}) yield $\si_m<\dfrac{1}{4}$ and the proof is complete. \qb
\end{Proof}
\begin{lem}
$\dlim_{\nu\ra+\infty}r_\nu=+\infty$.
\end{lem}
\begin{Proof}
Below, let us rewrite the relations that define the numbers $(r_\nu)$, $\nu=0,1,2,\ld$.
\begin{eqnarray}
r_1=r_0+\frac{c_2}{|\mi_{m_1(1)}|},  \label{eq3.9}
\end{eqnarray}
\begin{eqnarray}
r_2=r_1+\frac{c_2}{|\mi_{m_2}|},  \label{eq3.10}
\end{eqnarray}
\begin{eqnarray}
r_{\nu+1}=r_\nu+\frac{c_2}{|\mi_{m_{\nu+1}}|}, \ \ \nu=1,2,\ld\, ,  \label{eq3.11}
\end{eqnarray}
where $m_2:=m_1(m_1(1)+1)$, see subsection \ref{subsec3.3}. Equalities (\ref{eq3.9}), (\ref{eq3.10}),
(\ref{eq3.11}) imply
\begin{eqnarray}
r_\nu=r_0+c_2\cdot\sum^\nu_{k=1}\frac{1}{|\mi_{m_k}|} \ \ \text{for} \ \ \nu=1,2,\ld, \ \ \text{where} \ \
m_1=m_1(1).  \label{eq3.12}
\end{eqnarray}
By the definitions of $m_1(1)$, $m_2$ we have
\begin{eqnarray}
\sum^{m_1(1)}_{k=1}\frac{1}{|\mi_k|}\le
c_3\cdot\frac{1}{|\mi_1|}+\frac{1}{|\mi_{m_1(1)}|}<(c_3+1)\frac{1}{|\mi_1|},  \label{eq3.13}
\end{eqnarray}
\begin{eqnarray}
\sum^{m_2}_{k=m_1(1)+1}\frac{1}{|\mi_k|}\le
c_3\cdot\frac{1}{|\mi_{m_1(1)+1}|}+\frac{1}{|\mi_{m_2}|}<(c_3+1)\cdot\frac{1} {|\mi_{m_1(1)}|}.  \label{eq3.14}
\end{eqnarray}
Inductively, for every $\nu\ge2$ we get
\begin{eqnarray}
\sum^{m_{\nu+1}}_{k=m_\nu+1}\frac{1}{|\mi_k|}\le c_3\cdot\frac{1}{|\mi_{m_\nu+1}|}+\frac{1}{|\mi_{m_{\nu+1}}|}
<(c_3+1) \frac{1}{|\mi_{m_\nu}|}  \label{eq3.15}
\end{eqnarray}
because the sequence $(|\mi_n|)$ is strictly increasing. So by (\ref{eq3.13}), (\ref{eq3.14}) and (\ref{eq3.15})
we conclude that
\begin{eqnarray}
\ssum^{m_{\nu+1}}_{k=1}\dfrac{1}{|\mi_k|}<(c_3+1)\cdot\ssum^\nu_{k=0}\dfrac{1}{|\mi_{m_k}|}, \label{eq3.16}
\end{eqnarray}
where
$$
m_0:=1, \ \ m_1:=m_1(1).
$$
On the other hand, $\ssum^{+\infty}_{k=1}\dfrac{1}{|\mi_k|}=+\infty$ by our assumption. This fact and
(\ref{eq3.16}) give us
\begin{eqnarray}
\sum^{+\infty}_{k=0}\frac{1}{|\mi_{m_k}|}=+\infty.  \label{eq3.17}
\end{eqnarray}
Now by (\ref{eq3.12}) and (\ref{eq3.17}) we conclude that $\dlim_{\nu\ra+\infty}r_\nu=+\infty$ and the proof is
complete. \qb
\end{Proof}
\section{Construction and properties of the disks} \label{sec4}
\noindent

Fix the numbers $r_0,R_0,\thi_0,\thi_T,c_1,c_2,c_3,c_4$ which are defined in section \ref{sec2} and subsections
\ref{subsec3.1}, \ref{subsec3.2}, \ref{subsec3.3}. For the rest of this section we fix a subsequence $(\mi_n)$
of $(\la_n)$ satisfying the following:

1) $|\mi_n|, |\mi_{n+1}|-|\mi_n|>c_1$ for $n=1,2,\ld$ \bigskip

2) $\ssum^{+\infty}_{k=1}\dfrac{1}{|\mi_k|}=+\infty$.\bigskip

Finally, on the basis of the above, we consider the partition $\cp$ constructed in subsection \ref{subsec3.3}.
\subsection{Construction of the disks} \label{subsec4.1}
\noindent

Our goal in this subsection is to construct a certain family of pairwise disjoint disks, based on the previous
partition $\cp$ of the sector $S$. This family points out how one can use Runge's theorem to conclude the
Proposition \ref{prop2.1}. Let us describe, very briefly, the highlights of our argument. The main idea is to
assign to each point $w$ of the partition $\cp$ a suitable closed disk $B(w\mi(w),c_4)$ with center $w\mi(w)$
and radius $c_4$ (the radius will be the same for every member of the family of the disks), where $\mi(w)$ will
be chosen from the sequence $(\mi_n)$, so that on the one hand the disks $B(w\mi(w),c_4)$, $w\in\cp$ are
pairwise disjoint and on the other hand the union of the disks, $\cup_{w\in \cp}B(w\mi(w),c_4)$ "almost fills"
the sector $S$. It is evident that doing that, all the "good properties" of the partition established in the
previous section will
pass now to the family of the disks.\\

So, let us begin with the desired construction. We set
\[
\cb:=\{z\in\C/|z|\le c_4\}.
\]
Let $w\in\cp$ be a fixed point in $\cp$. By the definition of $\cp$ there exist unique
$r'\in\{r_0,r_1,\ld,r_{\nu_0}\}$, $m'\in\{1,m_1(1)+1,m_2+1,\ld,m_{\nu_0}+1\}$ such that $w\in\cp^{r',m'}_0$. By
definition, $\cp_0^{r',m'}=\phi_{r'}(\De_{m'})$. So there exists unique $n\in\{0,1,\ld,\nu_{m'}\}$ such that
$w=r'e^{2\pi i\thi^{m'}_n}$. Now there exist unique $k\in\N$, $k\ge1$ and $j\in\{0,1,\ld,m_1(m')-m'\}$ such that
$n=k(m_1(m')-m'+1)+j$, so we define $$\mi(w):=\mi_{m'+j}.$$ Thus we assign, in a unique way, a term of the
sequence $(\mi_n)$ to each one from the points of $\cp$. Finally we set $$\cb_w:=\cb+w\mi(w).$$ The desired
family of disks is the following:
$$\mathfrak{D}:=\{\cb\}\cup\{\cb_w:w\in\cp\}.$$
\subsection{Properties of the disks}\label{subsec4.2}
\begin{lem}\label{lem4.1}

We have $\cb\cap\cb_w=\emptyset$ for every $w\in\cp$.
\end{lem}
\begin{Proof}
$c_3=\dfrac{c_4}{r_0c_2}>\dfrac{c_4}{r_0}$, since $c_2\in(0,1)$. So $2c_3>2\dfrac{c_4}{r_0}$ and in view of
$c_1=4(c_3+1)>2c_3$ we get
\begin{eqnarray}
c_1>\frac{2c_4}{r_0}.  \label{eq4.1}
\end{eqnarray}
Take $w\in\cp$. The closed disks $\cb$, $\cb_w$ are centered at, 0, $w\mi(w)$ respectively and they have the
same radius $c_4$. Hence, we have to show that $|w\mi(w)|>2c_4$. Since $|w|\ge r_0$, it suffices to prove that
$|\mi(w)|>\dfrac{2c_4}{r_0}$. Observe now that, by the definition of $\mi(w)$ in the previous subsection,
\begin{eqnarray}
\mi(w)=\mi_n  \label{eq4.2}
\end{eqnarray}
for some positive integer $n\in\N$ and from the choice of $(\mi_n)$
\begin{eqnarray}
|\mi_n|>c_1 \ \ \text{for every} \ \ n\in\N.  \label{eq4.3}
\end{eqnarray}
Now, (\ref{eq4.1}), (\ref{eq4.2}) and (\ref{eq4.3}) imply $|\mi(w)|>2\dfrac{c_4}{r_0}$ and this finishes the
proof of the lemma. \qb
\end{Proof}
\begin{lem} \label{lem4.2}
Let $w_1,w_2\in\cp$ such that $|w_1|<|w_2|$. Then $\cb_{w_1}\cap\cb_{w_2}=\emptyset$.
\end{lem}
\begin{Proof}
We have
\[
m_0=1<m_1(1)+1,
\]
\[
m_2=m_1(m_1(1)+1)>m_1(1)
\]
and generally
\[
m_{\nu+1}=m_1(m_\nu+1)>m_\nu \ \ \text{for} \ \ \nu=1,2,\ld,\nu_0.
\]
Since $w_1,w_2\in\cp$, we have $w_1\in\cp_0^{r_{\nu_1},m_{\nu_1}+1}$, $w_2\in\cp_0^{r_{\nu_2},m_{\nu_2}+1}$ for
some $\nu_1,\nu_2\in\{0,1,\ld,\nu_0\}$ and so $|w_1|=r_{\nu_1}$, $|w_2|=r_{\nu_2}$. Our hypothesis
$|w_1|<|w_2|\Leftrightarrow r_{\nu_1}<r_{\nu_2}$ and the fact that the sequence $(r_\nu)$ is strictly increasing
gives us $\nu_1<\nu_2$. Thus, $m_{\nu_1}+1<m_{\nu_2}+1$, because the finite sequence $(m_\nu)$,
$\nu\in\{0,1,\ld,\nu_0\}$ is strictly increasing; recall that $m_0=1$, $m_1=m_1(1)$. By the definition of
$\mi(w)$ for $w\in\cp^{r',m'}_0\subset\cp$ we get $\mi(w)=\mi_{m'+j}$ for some $j\in\{0,1,\ld,m_1(m')-m'\}$, so
$|\mi_{m'}|\le|\mi(w)|\le|\mi_{m_1(m')}|$, since the sequence $(|\mi_n|)$ is strictly increasing. The fact that
$w_1\in\cp^{r_{\nu_1},m_{\nu_1}+1}_0$ implies
\begin{align*}
|\mi_{m_{\nu_1}+1}|&\le|\mi(w_1)|\le|\mi_{m_1(m_{\nu_1}+1)}| \\
&=|\mi_{m_{\nu_1+1}}|<|\mi_{m_{\nu_1+1}+1}|\le|\mi_{m_{\nu_2}+1}|,
\end{align*}
since $\nu_1+1\le \nu_2$ and the sequence $(|\mi_n|)$ is strictly increasing (\ref{eq4.1}). On the other hand we
have $w_2\in\cp_0^{r_{\nu_2},m_{\nu_2+1}}$, so
\[
|\mi_{m_{\nu_2}+1}|\le|\mi(w_2)|\le|\mi_{m_{\nu_2+1}}|.
\]
Hence, the last two inequalities above give
\[
|\mi(w_1)|<|\mi(w_2)|,
\]
which in turn implies
\begin{eqnarray}
|w_2\mi(w_2)|>|w_1\mi(w_1)|.  \label{eq4.4}
\end{eqnarray}
By (\ref{eq4.4}) and the hypothesis we get
\begin{align*}
|w_2\mi(w_2)-w_1\mi(w_1)|&\ge\big||w_2\mi(w_2)|-|w_1\mi(w_1)|\big|\\
&=|w_2\mi(w_2)|-|w_1\mi(w_1)|>|w_1|\,|\mi(w_2)|-|w_1|\,|\mi(w_1)| \\
&\ge r_0(|\mi(w_2)|-|\mi(w_1)|)>r_0c_1>2c_4,
\end{align*}
where the last inequality in the right hand side above follows from $c_1>\dfrac{2c_4}{r_0}$, which is already
established in Lemma \ref{lem4.1}. This shows that $\cb_{w_1}\cap\cb_{w_2}=\emptyset$. \qb
\end{Proof}
\begin{lem}\label{lem4.3}
Let $w_1,w_2\in\cp$ such that $w_1\neq w_2$ and $|w_1|=|w_2|$. Then $\cb_{w_1}\cap\cb_{w_2}=\emptyset$.
\end{lem}
\begin{Proof}
We distinguish two cases:\smallskip

(i) $|\mi(w_1)|<|\mi(w_2)|$. \smallskip

In this case, by our hypothesis, we have
\begin{align*}
|w_2\mi(w_2)-w_1\mi(w_1)|&\ge\big||w_2\mi(w_2)|-|w_1\mi(w_1)|\big|\\
&=|w_1|
\cdot(|\mi(w_2)|-|\mi(w_1)|)\ge r_0\cdot c_1>2c_4.
\end{align*}

Therefore $\cb_{w_1}\cap\cb_{w_2}=\emptyset$.\smallskip

(ii) $|\mi(w_1)|=|\mi(w_2)|$. \smallskip

Since $w_1,w_2\in\cp$ it follows that $w_1\in\cp_0^{r_{\nu_1},m_{\nu_1}+1}$,
$w_2\in\cp_0^{r_{\nu_2},m_{\nu_2}+1}$ for some $\nu_1,\nu_2\in\{0,1,\ld,\nu_0\}$. By the equalities
$|w_1|=r_{\nu_1}$, $|w_2|=r_{\nu_2}$ and the hypothesis $|w_1|=|w_2|$ we conclude that $r_{\nu_1}=r_{\nu_2}$,
which in turn implies $\nu_1=\nu_2$, since the sequence $(r_\nu)$ is strictly increasing. Setting
$\nu_1=\nu_2=\nu'$ we get $w_1,w_2\in\cp_0^{r_{\nu'},m_{\nu'}+1}$ for some $\nu'\in\{0,1,\ld,\nu_0\}$, that is
$w_1,w_2$ belong to the same partition of zero order. For simplicity we write $m_{\nu'}+1=m'$. We also set
$r_{\nu'}=r'$. So, $w_1,w_2\in\cp^{r',m'}_0$ and the definition of the set $\cp^{r',m'}_0$ gives us $w_1=r'\cdot
e^{2\pi i\thi^{(m')}_{n_1}}$, $w_2=r'\cdot e^{2\pi i\thi^{(m')}_{n_2}}$ for some
$n_1,n_2\in\{0,1,\ld,\nu_{m'}\}$, $n_1\neq n_2$, since $w_1\neq w_2$. Without loss of generality suppose that
$n_1<n_2$. Now, there exists a unique pair $(k_1,j_1)$, where $k_1\in\N$, $j_1\in\{0,1,\ld,m_1(m')-m'\}$ and a
unique pair $(k_2,j_2)$ where $k_2\in\N$ and $j_2\in\{0,1,\ld,m_1(m')-m'\}$ such that
\begin{eqnarray}
n_1=k_1(m_1(m')-m'+1)+j_1 \label{eq4.5}
\end{eqnarray}
and
\begin{eqnarray}
n_2=k_2(m_1(m')-m'+1)+j_2.  \label{eq4.6}
\end{eqnarray}
By definition, $\mi(w_1)=\mi_{m'+j_1}$ and $\mi(w_2)=\mi_{m'+j_2}$ and our hypothesis implies
\[
|\mi(w_1)|=|\mi(w_2)|\Leftrightarrow\mi(w_1)=\mi(w_2).
\]
So we have $j_1=j_2=j_0$ and
\[
\thi^{(m')}_{n_1}=\thi^{(m')}_{j_0}+k_1\si_{m'} ,
\]
\[
\thi^{(m')}_{n_2}=\thi^{(m')}_{j_0}+k_2\si_{m'}.
\]
Thus
\begin{eqnarray}
\thi^{(m')}_{n_2}-\thi^{(m')}_{n_1}=(k_2-k_1)\si_{m'}.  \label{eq4.7}
\end{eqnarray}
By (\ref{eq4.5}), (\ref{eq4.6}) and the fact that $n_1<n_2$ and $j_1=j_2$ we have $k_1<k_2\Rightarrow k_2\ge
k_1+1$. So, in view of (\ref{eq4.7}) we arrive at
\begin{eqnarray}
\thi^{(m')}_{n_2}-\thi^{(m')}_{n_1}\ge\si_{m'}>0.  \label{eq4.8}
\end{eqnarray}
The previous imply the following bound.
\begin{align}
|w_2\mi(w_2)-w_1\mi(w_1)|&=|\mi(w_1)|\cdot|w_1-w_2|\ge\mi_{m'}|\cdot|w_1-w_2| \nonumber\\
&=|\mi_{m'}|\cdot|r'\cdot e^{2\pi i\thi^{(m')}_{n_2}}-r'e^{2\pi i\thi^{(m')}_{n_1}}| \nonumber \\
&=r'|\mi_{m'}|\cdot|e^{2\pi i\thi^{(m')}_{n_2}}-e^{2\pi i\thi^{(m')}_{n_1}}| \nonumber \\
&=r'|\mi_{m'}|\cdot 2\sin(\pi(\thi^{(m')}_{n_2}-\thi^{(m')}_{n_1})) \nonumber\\
&\ge r_0\cdot|\mi_{m'}|\cdot2\sin(\pi(\thi^{(m')}_{n_2}-\thi^{(m')}_{n_1})).  \label{eq4.9}
\end{align}
Now, consider Jordan's inequality $$\sin x>\dfrac{2}{\pi}x, \,\,\,\, x\in\Big(0,\dfrac{\pi}{2}\Big).$$ We have
\[
0<\thi^{(m')}_{n_2}-\thi^{(m')}_{n_1}\leq \frac{1}{4}\Rightarrow0<\pi(\thi^{(m')}_{n_2}-
\thi^{(m')}_{n_1})<\frac{\pi}{4}.
\]
So, applying Jordan's inequality for
\[
x=\pi(\thi^{(m')}_{n_2}-\thi^{(m')}_{n_1})
\]
we get
\begin{eqnarray}
\sin(\pi(\thi^{(m')}_{n_2}-\thi^{(m')}_{n_1}))
>2(\thi^{(m')}_{n_2}-\thi^{(m')}_{n_1}).  \label{eq4.10}
\end{eqnarray}
By (\ref{eq4.8}), (\ref{eq4.9}) and (\ref{eq4.10}) it follows that
\begin{eqnarray}
|w_2\mi(w_2)-w_1\mi(w_1)|>4r_0|\mi_{m'}|\cdot\si_{m'}.  \label{eq4.11}
\end{eqnarray}
The definition of the number $\si_{m'}$ and relation (\ref{eq3.6}) of Lemma \ref{lemma3.4.1} yield
\[
\si_{m'}=c_2\cdot\sum^{m_1(m')}_{k=m'}\frac{1}{|\mi_k|}.
\]
By this fact, inequality (\ref{eq4.11}) and the definition of the number $m_1(m')$ we get
\begin{align}
|w_2\mi(w_2)-w_1\mi(w_1)|&>4r_0|\mi_{m'}|\cdot c_2\sum^{m_1(m')}_{k=m'}
\frac{1}{|\mi_k|} \nonumber \\
&>4r_0|\mi_{m'}|\cdot c_2\frac{c_3}{|\mi_{m'}|}=4r_0c_2c_3.  \label{eq4.12}
\end{align}
Recall that $c_3=\dfrac{c_4}{r_0c_2}$. So
\[
4r_0c_2c_3=4r_0c_2\cdot\frac{c_4}{r_0c_2}=4c_4>2c_4.
\]
The last bound along with (\ref{eq4.12}) give $\cb_{w_1}\cap\cb_{w_2}=\emptyset$ and the proof of the lemma is
complete. \qb
\end{Proof}

By Lemmas \ref{lem4.1}, \ref{lem4.2}, \ref{lem4.3} we conclude the following
\setcounter{cor}{3}
\begin{cor}\label{cor4.4}
The family $\mathfrak{D}:=\{\cb\}\cup\{\cb_w:w\in\cp\}$ consists of pairwise disjoint disks.
\end{cor}
\section{Proof of Lemma \ref{lem2.3}}\label{sec5}
\noindent

Let $j_1,s_1,k_1\in\N$ be fixed. Our aim is to prove that the set
$\bigcup\limits^{+\infty}_{m=1}E(m,j_1,s_1,k_1)$ is dense in $H(\mathbb{C})$. For simplicity we write
$p_{j_1}=p$. Fix $g\in H(\mathbb{C})$, a compact set $C\subseteq\C$ and $\e_0>0$. We seek $f\in H(\mathbb{C})$
and a positive integer $m_1$ such that
\begin{equation}\label{sec5eq1}
f\in E(m_1,j_1,s_1,k_1)
\end{equation}
and
\begin{equation} \label{sec5eq2}
\sup_{z\in C}|f(z)-g(z)|<\e_0.
\end{equation}
Fix $R_1>0$ sufficiently large so that
$$
C\cup\{z\in\C|\,|z|\le k_1\}\subset\{z\in\C|\,|z|\le R_1\}
$$
and then choose $0<\de_0<1$ such that
\begin{equation} \label{sec5eq3}
\textrm{if} \,\,\, |z|\le R_1 \,\,\, \textrm{and} \,\,\, |z-w|<\de_0, \ \ w\in\C, \,\,\, \textrm{then} \,\,\,
|p(z)-p(w)|<\frac{1}{2s_1}.
\end{equation}
Define $$\cb:=\{z\in\C|\,|z|\le R_1+\de_0\},$$
\[
c_4:=R_1+\de_0, \quad c_2:=\frac{\de_0}{2(2R_0\pi+1)},
\]
\[
c_3=\frac{c_4}{r_0c_2}=\frac{R_1+\de_0}{r_0\dfrac{\de_0}{2(2R_0\pi+1)}}=
\frac{2(R_1+\de_0)(2R_0\pi+1)}{r_0\de_0},
\]
\[
c_1=4(c_3+1)=4\cdot\bigg(\frac{2(R_1+\de_0)(2R_0\pi+1)}{r_0\de_0}+1\bigg).
\]
After the definition of the above numbers we choose a subsequence $(\mi_n)$ of $(\la_n)$ such that\smallskip

(i) $|\mi_n|>c_1, \,\,\, |\mi_{n+1}|-|\mi_n|>c_1$ for $n=1,2,\ld$ and \bigskip

(ii) $\ssum^{+\infty}_{n=1}\dfrac{1}{|\mi_n|}=+\infty$. \bigskip

On the basis of the fixed numbers $r_0,R_0,\thi_0,\thi_T,c_1,c_2,c_3,c_4$ and the choice of the sequence
$(\mi_n)$ we define the set $L$ as follows:
$$L:=\cb\cup \left( \bigcup_{w\in \mathcal{P}}\cb_w \right) ,$$
where the partition $\mathcal{P}$ and the discs $\cb_w$, $w\in \cp$ are constructed in Sections 3 and 4
respectively. By Corollary \ref{cor4.4}, the family $\mathfrak{D}$ consists of pairwise disjoint disks.
Therefore the compact set $L$ has connected complement. This property is needed in order to apply Mergelyan's
theorem. We now define the function $h$ on the compact set $L$, $h:L\ra\C$ by
\[
h(z)=\left\{\begin{array}{cc}
              g(z), & z\in \cb \\
              p(z-w\lambda (w)), & z\in \cb_w ,w\in \cp.
            \end{array}\right.
\]
By Mergelyan's theorem \cite{12} there exists an entire function $f$ (in fact a polynomial) such that
\begin{eqnarray}
\sup_{z\in L}|f(z)-h(z)|<\min\bigg\{\frac{1}{2s_1},\e_0\bigg\}.  \label{sec5eq4}
\end{eqnarray}
The definition of $h$ and (\ref{sec5eq4}) give
\[
\sup_{z\in C}|f(z)-g(z)|\le\sup_{z\in\cb}|f(z)-g(z)|=\sup_{z\in L}|f(z)-h(z)|<\e_0,
\]
which implies the desired inequality (\ref{sec5eq2}). It remains to show (\ref{sec5eq1}).

Let $a\in S$. Then $a=re^{2\pi i\thi}$ for some $r\in[r_0,R_0]$ and $\thi\in[\thi_0,\thi_T]$. There exists a
unique $n_0\in\{0,1,\ld,\nu_0-1\}$ such that either $r_{n_0}\le r<r_{n_0+1}$ or $r_{\nu_0}\le r\le R_0$. We set
$$r_1:=r_{n_0}, \,\,\, r_2:=r_{n_0+1} \,\,\, \textrm{if}\,\,\, r_{n_0}\le r<r_{n_0+1}$$
and
$$r_1:=r_{\nu_0},\,\,\, r_2:=R_0 \,\,\,\textrm{if}\,\,\, r_{\nu_0}\le r\le R_0.$$

By the construction of the partition $\cp$ there exists a unique $m'\in\N$ such that $\cp^{r_1,m'}_0\subset\cp$.
In addition, there exists unique $\rho\in\{0,1,\ld,\nu_{m'}-1\}$ such that
\[
\textrm{either}\,\,\, \thi^{(m')}_\rho\le\thi<\thi^{(m')}_{\rho+1} \ \ \text{or} \ \
\thi^{(m')}_{\nu_{m'}}\le\thi\le\thi_T.
\]
Define now  $$\thi_1:=\thi^{(m')}_\rho, \thi_2:=\thi^{(m')}_{\rho+1} \,\,\, \textrm{if}\,\,\,
\thi^{(m')}_\rho\le\thi<\thi^{(m')}_{\rho+1} $$ and $$\thi_1:=\thi^{(m')}_{\nu_{m'}}, \thi_2:=\thi_T \,\,\,
\textrm{if} \,\,\, \thi^{(m')}_{\nu_{m'}}\le\thi\le\thi_T$$ and then set $$w_0:=r_1\cdot e^{2\pi
i\thi_1}\in\cp_0^{r_1,m'}.$$

We shall prove now that for every $z\in\C$ with $|z|\le R_1$, $z+a\mi(w_0)\in\cb_{w_0}$. Recall that
$\cb_{w_0}:=\cb+w_0\mi(w_0)=\overline{D}(w_0\mi(w_0),R_1+\de_0)$. It suffices to prove that
\begin{eqnarray}
|(z+a\mi(w_0))-w_0\mi(w_0)|<R_1+\de_0 \ \ \text{for} \ \ |z|\le R_1.  \label{sec5eq5}
\end{eqnarray}
For $|z|\le R_1$ we have,
\begin{align}
|(z+a\mi(w_0))-w_0\mi(w_0)|&\le R_1+|\mi(w_0)|\,|a-w_0|\nonumber \\
&=R_1+|\mi(w_0)|\cdot| r\cdot e^{2\pi i\thi}-r_1e^{2\pi i\thi_1}|.  \label{sec5eq6}
\end{align}
By (\ref{sec5eq5}) and (\ref{sec5eq6}) it suffices to prove
\begin{eqnarray}
|\mi(w_0)|\cdot|re^{2\pi i\thi}-r_1e^{2\pi i\thi_1}|<\de_0.  \label{sec5eq7}
\end{eqnarray}
We now have
\begin{align*}
|re^{2\pi i\thi}-r_1e^{2\pi i\thi_1}|&=|re^{2\pi i\thi}-r_1e^{2\pi i\thi}+
r_1e^{2\pi i\thi}-r_1e^{2\pi i\thi_1}| \nonumber\\
&\le|re^{2\pi i\thi}-r_1e^{2\pi i\thi}|+|r_1e^{2\pi i\thi}-r_1e^{2\pi i\thi_1}| \nonumber\\
&\le|r-r_1|+r_1|e^{2\pi i\thi}-e^{2\pi i\thi_1}| \nonumber\\
&\le|r_1-r_2|+R_0\cdot2\sin(\pi(\thi_1-\thi)) \nonumber \\
&\le(r_2-r_1)+R_02\sin(\pi(\thi_2-\thi_1)) \nonumber\\
&<(r_2-r_1)+2R_0\pi(\thi_2-\thi_1) \nonumber \\
&\le2 R_0\pi\cdot\frac{c_2}{|\mi(w_0)|}+\frac{c_2}{|\mi(w_0)|} \nonumber \\
&=(2R_0\pi+1)c_2\frac{1}{|\mi(w_0)|} \nonumber \\
&=(2R_0\pi+1)\cdot\frac{\de_0}{2(2R_0\pi+1)}\cdot\frac{1}{|\mi(w_0)|}=\frac{\de_0} {2|\mi(w_0)|}
\end{align*}
which implies (\ref{sec5eq7}). So we proved that for every $z\in\C$, $|z|\le R_1$
\begin{eqnarray}
z+a\mi(w_0)\in\cb_{w_0}.  \label{sec5eq8}
\end{eqnarray}
By the definition of $h$ and (\ref{sec5eq8}) we have that for every $z\in\C$ with $|z|\le R_1$
\begin{eqnarray}
|f(z+a\mi(w_0))-p(z+\mi(w_0)(re^{2\pi i\thi}-r_1e^{2\pi i\thi_1})|<\frac{1}{2s_1}.  \label{sec5eq9}
\end{eqnarray}
Take any $z\in\C$ with $|z|\le R_1$. By (\ref{sec5eq3}) and (\ref{sec5eq7})
\begin{eqnarray}
|p(z+\mi(w_0)(re^{2\pi i\thi}-r_1e^{2\pi i\thi_1}))-p(z)|<\frac{1}{2s_1}  \label{sec5eq10}
\end{eqnarray}
and the triangle inequality gives
\begin{align}
|f(z+a\mi(w_0))-p(z)|\le&|f(z+a\mi(w_0))-p(z+\mi(w_0)(re^{2\pi i\thi}-r_1
e^{2\pi i\thi_1}))| \nonumber\\
&+|p(z+\mi(w_0)(re^{2\pi i\thi}-r_1e^{2\pi i\thi_1}))-p(z)|.  \label{sec5eq11}
\end{align}
Using (\ref{sec5eq9}), (\ref{sec5eq10}), (\ref{sec5eq11}) we arrive at
\[
|f(z+a\mi(w_0))-p(z)|<\frac{1}{s_1}
\]
and since $k_1\le R_1$ it readily follows that
\begin{eqnarray}
\sup_{|z|\le k_1}|f(z+a\mi(w_0))-p(z)|<\frac{1}{s_1}.  \label{sec5eq12}
\end{eqnarray}
Set
\[
m_1:=\max\{n\in\N|\la_n=\mi(w), \ \ \text{for some} \ \ w\in\cp\}
\]
and observe that the definition of $m_1$ is independent from $a\in S$. Thus, by the previous we conclude that
for every $a\in S$ there exists some $n\in\N$, $n\le m_1$ such that
\[
\sup_{|z|\le k_1}|f(z+a\la_n)-p(z)|<\frac{1}{s_1},
\]
where $f\in H(\mathbb{C})$, since $f$ is a polynomial. This completes the proof of the lemma. \qb
\section{Examples of sequences $\bbb{\La:=(\la_n)}$ satisfying condition $\bbb{\Si}$}
\noindent

By the remark in \cite{7} we have a sample of first examples satisfying condition $(\Si)$:
\[
\la_n=n, \ \ \la_n=n(\log n)^p \ \  \text{for} \ \ p\le1,  \ \ \la_n=n\log n\log\log n.
\]
In all the above examples we also have $\Big|\dfrac{\la_{n+1}}{\la_n}\Big|\ra1$ as $n\ra+\infty$. However, for
sequences $(\la_n)$, such that $\la_n\ra\infty$ and $\Big|\dfrac{\la_{n+1}}{\la_n}\Big|\ra1$ we have that the
conclusion of Theorem \ref{thm2.1} holds by the main result in \cite{Tsi}. It is our aim to show that there
exist sequences $(\la_n)$, such that: $\la_n\ra\infty$, $(\la_n)$ satisfies condition $(\Si)$ and the ratio
$\Big|\dfrac{\la_{n+1}}{\la_n}\Big|$ does not tend to $1$.

Let us see things more specifically. Consider a sequence  $\La=(\la_n)$ of non-zero complex numbers and define
\[
\cb(\La):=\bigg\{a\in[0,+\infty]|\exists\,(\mi_n)\subset\La\ \ \text{with} \ \
a=\underset{n}{\lim\sup}\bigg|\frac{\mi_{n+1}}{\mi_n}\bigg|\bigg\},
\]
$$i(\La):=\inf\cb(\La) .$$
Clearly $$i(\La)\in[0,+\infty]$$ and $$\textrm{if}\,\,\, \la_n\ra\infty \,\,\,\textrm{then} \,\,\,
\cb(\La)\subset[1,+\infty] \,\,\, \textrm{and} \,\,\, i(\La)\in[1,+\infty].$$

We say that a sequence of non-zero complex numbers $\Lambda =(\lambda_n)$ satisfies condition $(\Si ')$ if
$i(\Lambda )=1$. In \cite{Tsi} we established the following result.\medskip

\textit{If $\La=(\la_n)$ is a sequence of non-zero complex numbers such that $\la_n\ra\infty$ and $\Lambda$
satisfies condition $(\Si ')$, then the conclusion of Theorem \ref{thm2.1} holds}.\medskip

In view of the above result the following question arises naturally.
\begin{que} \label{Q2}
Let $\La=(\la_n)$ be a sequence of non-zero complex numbers such that $\la_n\ra\infty$ and $i(\La)>1$. Does the
conclusion of Theorem \ref{thm2.1} hold?
\end{que}
It is quite surprising, at least to us, that the answer to the above question is sometimes yes and sometimes
no!. In what follows we shall exhibit examples of sequences admitting a positive answer to this question.
Results going to the opposite direction are established in \cite{8}. In particular, the main result in \cite{8}
is the following: if $\La=(\la_n)$ is a sequence of non-zero complex numbers such that
$\liminf_{n}\frac{|\la_{n+1}|}{|\la_n|}>2$ (hence $i(\La)>2$) then
$\bigcap\limits_{a\in\C\sm\{0\}}HC(\{T_{\la_na}\})=\emptyset$.

Below we construct specific examples of sequences $\La=(\la_n)$ such that $\la_n\to \infty$, $i(\La)=M$ for any
fixed positive number $M>1$ and $\La$ satisfies $(\Si)$. By this result we complete our goals in this paper that
are the following three.

\begin{itemize}
\item Firstly, we give affirmative reply to Question 1 of \cite{7}. \item Secondly, for certain sequences, we
also give a positive answer to Question \ref{Q2}. \item Thirdly, we exhibit a variety of examples of sequences
$\Lambda =(\lambda_n)$ of non-zero complex numbers with $\lambda_n\to \infty$ such that $\Lambda$ satisfies
condition $(\Si )$ and it does not satisfy condition $(\Si ')$.
\end{itemize}

The above discussion shows that the problem of deciding whether a sequence $\La=(\lambda_n)$, such that
$\lambda_n\ra\infty$ and $i(\La)=M$ for some $M>1$ satisfies the conclusion of Theorem \ref{thm2.1} is quite
delicate and needs further study.
\begin{prop}\label{prop6.1}
For every $M>1$ there exists a sequence $\La=(\la_n)$ such that $\la_n\ra\infty$, $i(\La)=M$ and condition
$(\Si)$ holds for $\La$. Thus, for every $M>1$ there exists a sequence of non-zero complex numbers $\La=(\la_n)$
such that $\la_n\ra\infty$ as $n\ra+\infty$, $i(\La)=M$ and $\bigcap\limits_{a\in\C\sm\{0\}}HC(\{T_{\la_na}\})$
is a $G_\delta$ and dense subset of $H(\mathbb{C})$.
\end{prop}
\begin{Proof}
Fix a positive number $M_0>1$. We shall construct a sequence of non-zero complex numbers $\La=(\la_n)$ such that
$\la_n\ra\infty$, $i(\La)=M_0$ and condition $(\Si)$ holds for $\La$. The sequence $\La$ will be a strictly
increasing sequence of positive numbers such that $\la_n\ra+\infty$ as $n\ra+\infty$.

We construct inductively a countable family $\{\mathfrak{D}_n\}$, $n=1,2,\ld$ of sets $\mathfrak{D}_n \subset
[1,+\infty )$ according to the following rules.
\begin{enumerate}
\item[(i)] $\mathfrak{D}_1=\{ 1\}$.
\item[(ii)] $\mathfrak{D}_n=\{a_n+\nu|\nu=0,1,\ld,\;([a_n]+1)!\}, \,\,\,n=1,2,\ldots$.
\item[(iii)] $\min\mathfrak{D}_{n+1}=M_0\cdot\max\mathfrak{D_n}$ for each $n=1,2,\ldots $,
\end{enumerate}
where $a_n=\min\mathfrak{D_n}$ and $[x]$ denotes the integer part of the real number $x$ as usual. Observe that
every $n$, $m\in\N$, $n\neq m$, $\mathfrak{D}_n\cap\mathfrak{D_m}=\emptyset$. Set
$$\widetilde{\La}=\bigcup\limits^{+\infty}_{n=1}\mathfrak{D}_n.$$ We define the sequence $\La=(\la_n)$ to be the enumeration of
$\widetilde{\La}$ by the natural order.

It is obvious that $\la_n\neq0$ $\fa\;n\in\N$, $\dis\lim_{n\ra+\infty}\la_n=+\infty$, and $(\la_n)$ is a
strictly increasing sequence of positive numbers. We prove now the following
\end{Proof}
{\bf Claim 1:} For every subsequence $\mi=(\mi_n)$ of $\La$ we have $\underset{n\ra+\infty}{\lim\sup}\dfrac{\mi_{n+1}}{\mi_n}\ge M_0$.
\begin{Proof}
Firstly we prove that for every natural number $m\in\N$, there exists some $N\in\N$, $N\ge m$ such that
\[
\frac{\mi_{N+1}}{\mi_N}\ge M_0.
\]
So, take any $m\in\N$ and let $m_1$ be the unique natural number such chat $\mi_m\in\mathfrak{D}_{m_1}$. Setting
$A_{m_1}:=\{n\in\N|\mi_n\in\mathfrak{D}_{m_1}\}$, it is obvious that $A_{m_1}\neq\emptyset$, since $m\in
A_{m_1}$. We set $m_2:=\max A_{m_1}$. Then $\mi_{m_2+1}\notin\mathfrak{D}_{m_1}$ and so
$\mi_{m_2+1}\ge\min\mathfrak{D}_{m_1+1}$. We have $\mi_{m_2}\le\max\mathfrak{D}_{m_1}$, thus
\[
\frac{\mi_{m_2+1}}{\mi_{m_2}}\ge\frac{\min\mathfrak{D}_{m_1+1}}{\max\mathfrak{D}_{m_1}}
=M_0 \ \ \text{and} \ \ m_2\ge m_1.
\]
So we proved that for every $m\in\N$, there exists some $N\ge m$ such that $\dfrac{\mi_{N+1}}{\mi_N}\ge M_0$. We
incorporate the last fact into an inductive argument and obtain the following. For $m=1$ there exists
$k_1\in\N$, $k_1\ge1$ such that $\dfrac{\mi_{k_1+1}}{\mi_{k_1}}\ge M_0$. For $m=k_1+1$, there exists some
$k_2\ge k_1+1$ (especially $k_2>k_1$) such that $\dfrac{\mi_{k_2+1}}{\mi_{k_2}}\ge M_0$. Suppose that for some
$\nu\in\N$ we have found some $k_\nu\in\N$ such that $\dfrac{\mi_{k_\nu+1}}{\mi_{k_\nu}}\ge M_0$. Then for
$m=k_\nu+1$ there exists some $k_{\nu+1}\ge k_\nu+1$ (especially $k_{\nu+1}>k_\nu$) such that
$\dfrac{\mi_{k_{\nu+1}}+1}{\mi_{k_{\nu+1}}}\ge M_0$. Therefore we obtain a subsequence $(\mi_{k_\nu})$,
$\nu=1,2,\ld$ of $(\mi_n)$ such that $k_{\nu+1}>k_\nu$ for each $\nu=1,2,\ld$ and
$\dfrac{\mi_{k_\nu+1}}{\mi_{k_\nu}}\ge M_0$. This gives
$\underset{\nu\ra+\infty}{\lim\sup}\dfrac{\mi_{k_\nu+1}}{\mi_{k_\nu}}\ge M_0,$ which in turn implies
\[
\underset{n\ra+\infty}{\lim\sup}\frac{\mi_{n+1}}{\mi_n}\ge M_0.
\]
This completes the proof of Claim 1.
\end{Proof}
{\bf Claim 2:} $\underset{n\ra+\infty}{\lim\sup}\dfrac{\la_{n+1}}{\la_n}=M_0$.
\begin{Proof}
Let $n\in\N$. If $\la_n,\la_{n+1}\in\mathfrak{D}_m$ for some positive integer $m$, then by the construction of
$\mathfrak{D}_m$ we have
\begin{equation} \label{sec6eq1}
\la_{n+1}=\la_n+1\Rightarrow\dfrac{\la_{n+1}}{\la_n}=1+\dfrac{1}{\la_n}.
\end{equation}
If there is no $m\in\N$ such that $\la_n$, $\la_{n+1}\in\mathfrak{D}_m$, then this can happen only if
$\la_n=\max\mathfrak{D}_m$ and $\la_{n+1}=\min\mathfrak{D}_{m+1}$ for some $m\in\N$, hence
\begin{equation} \label{sec6eq2}
\dfrac{\la_{n+1}}{\la_n}=M_0.
\end{equation}
By (\ref{sec6eq1}), (\ref{sec6eq2}) and since $\la_n\to +\infty$ the conclusion follows. This completes the
proof of Claim 2.
\end{Proof}
Claims 1 and 2 imply that $i(\La)=M_0$.\\

{\bf Claim 3:} The sequence $\La$ satisfies condition $(\Si)$.
\begin{Proof}
Fix some natural number $N_0\ge2$. We will show that there exists a subsequence $(\mi_n)$ of $\La$ such that

(i) $\mi_{n+1}-\mi_n>N_0$ for every $n=1,2,\ld$ and

(ii) $\ssum^{+\infty}_{n=1}\dfrac{1}{\mi_n}=+\infty$.

Recall that $a_n=\min\mathfrak{D}_n>1$ for every $n\ge2$. Since
\[
\dis\lim_{n\ra+\infty}\Big(1+\frac{1}{2}+\cdots+\frac{1}{n}-\log n\Big)=\ga ,
\]
where $\ga\simeq0,57722156649\ldots$ is the Euler constant, there exists some natural number $n_0\in\N$ such
that
\begin{align*}
-\frac{1}{2}<\sum^n_{k=1}\frac{1}{k}-\log n-\ga<\frac{1}{2} \ \ \text{for} \ \ n\ge n_0>2.
\end{align*}
Let some $m,n\in\N$, $m>n\ge n_0$. Then we have
\begin{align}
\frac{1}{n+1}+\frac{1}{n+2}+\cdots+\frac{1}{m}=&\sum^m_{k=n+1}\frac{1}{k}
=\sum^m_{k=1}\frac{1}{k}-\sum^n_{n=1}\frac{1}{k}\nonumber \\
=&\bigg(\sum^m_{k=1}\frac{1}{k}-\log m-\ga\bigg)-\bigg(\sum^n_{k=1}\frac{1}{k}-\log n-\ga\bigg)  \nonumber\\
&+\log m-\log n>\log\frac{m}{n}-1 \nonumber \\
=&\log\frac{m}{n}+\log e^{-1}=\log\bigg(\frac{m}{n}\cdot e^{-1}\bigg)\nonumber\\
=&\log\bigg(\frac{m}{ne}\bigg).  \label{sec6eq3}
\end{align}
It is easy to show that  $ a_n>n$ for $n\ge 2$. Set $n_1:=\max\{n_0,N_0\}+2$. Let now some $n\in\N$ with $n\ge
n_1$. Recall that
\begin{align*}
\mathfrak{D}_n&=\{a_n,a_n+1,\ld,a_n+([a_n]+1)!\} \\
&=\{a_n+j|j=0,1,\ld,\;([a_n]+1)!\}
\end{align*}
Setting $N_1:=N_0+1$ we obtain
\begin{align}
&\frac{1}{a_n}+\frac{1}{a_n+N_1}+\frac{1}{a_n+2N_1}+\cdots+\frac{1}
{a_n+\dfrac{([a_n]+1)!}{N_1}\cdot N_1}\nonumber \\
&>\frac{1}{N_1a_n}+\frac{1}{N_1a_n+N_1}+\frac{1}{N_1a_n+2N_1}+\cdots+
\frac{1}{N_1a_n+N_1\cdot\dfrac{([a_n]+1)!}{N_1}} \nonumber\\
&=\frac{1}{N_1}\cdot\sum^{\frac{([a_n]+1)!}{N_1}}_{j=0}\frac{1}{a_n+j}>
\frac{1}{N_1}\cdot\sum^{\frac{([a_n]+1)!}{N_1}}_{j=0}\frac{1}{([a_n]+1)+j}. \label{sec6eq4}
\end{align}
We write for simplicity $\nu=[a_n]+1$. So by (\ref{sec6eq3}), (\ref{sec6eq4}) we get

\begin{align} \label{sec6eq4*}
\sum^{\frac{\nu!}{N_1}}_{k=0}\frac{1}{a_n+kN_1}>\frac{1}{N_1}\cdot\log\bigg(
\frac{\nu+\frac{\nu!}{N_1}}{(\nu-1)e}\bigg)
>\frac{1}{N_1}\cdot
\log\bigg(\frac{(\nu-1)!}{N_1e}\bigg).
\end{align}
We will show that
\[
\frac{1}{N_1}\cdot\log\bigg(\frac{(\nu-1)!}{N_1e}\bigg)>\nu
\]
for $\nu$ big enough. It follows that
\begin{align*}
\frac{1}{N_1}\cdot\log\bigg(\frac{(\nu-1)!}{N_1e}\bigg)>\nu&\Leftrightarrow\log
\bigg(\frac{(\nu-1)!}{N_1e}\bigg)>N_1\nu \\
&\Leftrightarrow(\nu-1)!>N_1e\cdot e^{N_1\nu}=N_1\cdot e^{N_1\nu+1}.
\end{align*}
Let us consider the sequence $\ga_\nu=\dfrac{(\nu-1)!}{N_1e^{N_1\nu+1}}$. By the ratio criterion for $(\ga_\nu)$
we have
\[
\frac{\ga_{\nu+1}}{\ga_\nu}=\frac{\dfrac{\nu!}{N_1e^{N_1(\nu+1)+1}}}{\dfrac{(\nu-1)!}
{N_1e^{N_1\nu+1}}}=\frac{\nu!\cdot e^{N_1\nu+1}}{(\nu-1)!\cdot e^{N_1(\nu+1)+1}}=\frac{\nu}
{e^{N_1}}.
\]
So $\dis\lim_{\nu\ra+\infty}\bigg(\dfrac{\ga_{\nu+1}}{\ga_\nu}\bigg)=+\infty$ which implies that there exists
some $n_2\ge n_1$ such that $\ga_n>1$ for $n\ge n_2$ or equivalently
\begin{eqnarray} \label{sec6eq6}
\dfrac{1}{N_1}\cdot\log\bigg(\dfrac{(n-1)!}{N_1e}\bigg)>n,\,\,\, n\ge n_2.
\end{eqnarray}
Thus by (\ref{sec6eq4*}) and (\ref{sec6eq6}) we have:
\[
\sum^{\frac{\nu!}{N_1}}_{k=0}\frac{1}{a_n+kN_1}>[a_n]+1 \ \ \text{for} \ \ n\ge n_2.
\]
Now for $n\ge n_2$ define the set
\[
\mathfrak{D}'_n:=\bigg\{a_n,a_n+N_1,a_n+2N_1,\ld,a_n+\frac{([a_n]+1)!}{N_1}\cdot N_1\bigg\},
\]
and consider the union
\[
\mathfrak{D}':=\bigcup\limits_{n\ge n_2}\mathfrak{D}'_n.
\]
Let $(\mi_n)$ be the sequence we get when we enumerate $\mathfrak{D}'$ by its natural order. Clearly $(\mi_n)$
is a subsequence of $\La$ and satisfies the desired properties (i) and (ii). This completes the proof of Claim 3
and hence that of Proposition \ref{prop6.1} using Theorem \ref{thm2.1}.  \qb
\end{Proof}

\begin{cor}
There exists a sequence $\Lambda =(\lambda_n)$ of non-zero complex numbers with $\lambda_n\to \infty$ such that
$\Lambda$ satisfies condition $(\Si)$ and it does not satisfy condition $(\Si ')$.
\end{cor}
\begin{Proof}
Every sequence $\Lambda =(\lambda_n)$ of non-zero complex numbers with $\lambda_n\to \infty$ which satisfies the
conclusion of Proposition \ref{prop6.1}, clearly does not satisfy $(\Si ')$.  \qb
\end{Proof}
We point out that sequences of the form $(n^2)$, $(n^3)$, $(n^4) \ldots $, satisfy condition $(\Si ')$ but they
do not satisfy $(\Si )$. To complete the picture we observe that there are sequences with sufficiently slow
growth, such as $(n)$, $(\sqrt{n} )$, $(\log (n+1) )$, $(\log \log (n+1))$, that satisfy both conditions $(\Si
)$ and $(\Si ')$. Hence, neither $(\Si )$ nor $(\Si ')$ implies the other and, in addition, they have non-empty
intersection. This in turn shows that Theorem \ref{thm2.1} does not follow by the main result in \cite{Tsi} and
vice versa.\medskip

We close the paper with a question which kindly posed to us by the referee.
\begin{que}
If $(\lambda_n)$ is a sequence of non-zero complex numbers, $\lambda_n\to \infty$ such that $\bigcap_{a\in
\mathbb{C}\setminus \{ 0\}} HC(\{ T_{\la_na} \}) \neq \emptyset $ what can be said about the growth of common
hypercyclic entire functions for the collection of sequences $(T_{\la_na} )$, $a\in \mathbb{C}\setminus \{ 0\}$?
To answer such a question, we should specify the sequence $(\lambda_n)$. For instance, what happens when
$\lambda_n=n^2$, $\lambda_n=n\log (n+1)$, $\lambda_n=n^3$, etc.?
\end{que}

Results concerning permissible and optimal growth rates for hypercyclic entire functions with respect to the
translation operator $T_a$, $a\in \mathbb{C}\setminus \{ 0\}$, as well as similar results for differential
operators acting on various function spaces can be found in \cite{AA1}, \cite{AA2}, \cite{Arm}, \cite{BeiMul},
\cite{BerBo}, \cite{BerBoCo}, \cite{BBG}, \cite{ChaSha}, \cite{DS}, \cite{GMPR}, \cite{Duios}, \cite{Grosse1},
\cite{Grosse3}, \cite{Shkarin}.
\medskip

{\bf Acknowledgements}. I am grateful to George Costakis for his helpful comments and remarks and for all the
help he offered me concerning the presentation of this work. I would also like to thank the referee for useful
comments which led me to improve the presentation of the present paper.

%
%
%
%
%
%
%

\end{document}